\documentclass[12pt,reqno]{amsart}

\usepackage{amssymb,bbm,esint,units,url}
\usepackage{booktabs}
\usepackage[bookmarks=false,unicode,colorlinks,urlcolor=red,citecolor=red,linkcolor=blue]{hyperref}
\usepackage{aliascnt}
\usepackage{doi}
\usepackage{tikz}
\usepackage{pgfplots} 

\usepackage{color}

\usepackage[normalem]{ulem}
\usepackage[margin=1.25in]{geometry}

\newtheorem{thmx}{Theorem}

\newaliascnt{propx}{thmx}

\aliascntresetthe{propx}

\newaliascnt{corx}{thmx}
\newtheorem{corx}[corx]{Corollary}
\aliascntresetthe{corx}

\newaliascnt{conjx}{thmx}

\aliascntresetthe{conjx}

\newaliascnt{lemma}{theorem}
\newtheorem{lemma}[lemma]{Lemma}
\aliascntresetthe{lemma}

\newaliascnt{proposition}{theorem}
\newtheorem{proposition}[proposition]{Proposition}
\aliascntresetthe{proposition}

\newaliascnt{corollary}{theorem}

\aliascntresetthe{corollary}

\newaliascnt{conjecture}{theorem}
\newtheorem{conjecture}[conjecture]{Conjecture}
\aliascntresetthe{conjecture}

\newaliascnt{example}{theorem}

\aliascntresetthe{example}

\makeatletter
\def\tagform@#1{\maketag@@@{\ignorespaces#1\unskip\@@italiccorr}}
\let\orgtheequation\theequation
\def\theequation{(\orgtheequation)}
\makeatother
\def\equationautorefname~{}

\newcommand{\arxiv}[1]{%
 \href{http://front.math.ucdavis.edu/#1}{ArXiv:#1}}

\newcommand{\B}{{\mathbb B}}

\newcommand{\D}{{\mathbb D}}

\newcommand{\R}{{\mathbb R}}

\newcommand{\Rn}{{\mathbb R}^n}

\newcommand{\ds}{\displaystyle}

\newcommand{\fr}[2]{\frac{\ds #1}{\ds #2}}

\begin{document}

\title{From Steklov to Neumann and beyond, via Robin: the Szeg\H{o} way}
\author[]{Pedro Freitas and Richard S. Laugesen}
\address{Departamento de Matem\'atica, Instituto Superior T\'ecnico, Universidade de Lisboa, Av. Rovisco Pais 1,
P-1049-001 Lisboa, Portugal {\rm and}
Grupo de F\'isica M\'atematica, Faculdade de Ci\^encias, Universidade de Lisboa,
Campo Grande, Edif\'icio C6, P-1749-016 Lisboa, Portugal}
\email{psfreitas\@@fc.ul.pt}
\address{Department of Mathematics, University of Illinois, Urbana,
IL 61801, U.S.A.}
\email{Laugesen\@@illinois.edu}
\date{\today}


\begin{abstract}
The second eigenvalue of the Robin Laplacian is shown to be maximal for the disk among simply-connected planar domains of fixed area 
when the Robin parameter is scaled by perimeter in the form $\alpha/L(\Omega)$, and $\alpha$ lies between $-2\pi$ and $2\pi$. Corollaries include Szeg\H{o}'s sharp upper bound on the second eigenvalue of the Neumann Laplacian under area normalization, and Weinstock's inequality for the first nonzero Steklov eigenvalue for simply-connected domains of given perimeter.

The first Robin eigenvalue is maximal, under the same conditions, for the degenerate rectangle. When area normalization on the domain is changed to conformal mapping normalization and the Robin parameter is positive, the maximiser of the first eigenvalue changes back to the disk. 
\end{abstract}

\maketitle

\section{\bf Introduction}
\label{intro}

The eigenvalue problem for the Robin Laplacian on a domain $\Omega \subset \R^2$ with Lipschitz boundary is
\begin{equation}\label{robinproblem}
\begin{split}
- \Delta u & = \lambda u \ \quad \text{in $\Omega$,} \\
\frac{\partial u}{\partial\nu} + \alpha u & = 0 \qquad \text{on $\partial \Omega$,} 
\end{split}
\end{equation}
where $\alpha$ is a real parameter and $\nu$ is the outward unit normal. The corresponding eigenvalues, denoted $\lambda_{k}(\Omega;\alpha)$ for $k=1,2,
\dots$, are increasing and continuous as functions of the Robin parameter $\alpha$, and for each fixed $\alpha$ satisfy
\[
\lambda_1(\Omega;\alpha) < \lambda_2(\Omega;\alpha) \leq \lambda_3(\Omega;\alpha) \leq \dots \to \infty . 
\]
Isoperimetric eigenvalue inequalities in the literature typically assume an area normalization of the domain --- for instance see \cite{B86,D06,FL18a} and the survey \cite{BFK}, which includes many related results on Robin eigenvalues.

While this area normalization is natural for Dirichlet and Neumann problems, it provides only part of the story for Robin because the rescaling relation $t^2 \lambda(t\Omega;\alpha/t) = \lambda(\Omega;\alpha)$ shows that the area-normalized product $|\Omega| \lambda(\Omega;\alpha)$ is not scale-invariant. This observation has prompted us to look for natural, scale-invariant isoperimetric inequalities for eigenvalues of problem \eqref{robinproblem}. We claim the most natural formulation for planar domains is to keep the domain normalized by area while considering the Robin parameter scaled by the perimeter of the domain. The eigenvalues under consideration thus become
\[
\lambda_1 \! \left(\Omega;\fr{\alpha}{L(\Omega)}\right) < \lambda_2 \! \left(\Omega;\fr{\alpha}{L(\Omega)}\right) \leq \lambda_3 \! \left(\Omega;\fr{\alpha}{L(\Omega)}\right) \leq \dots \to \infty,
\]
where $L(\Omega)$ denotes the length of the boundary $\partial \Omega$.
Under this new scaling, the behavior of eigenvalues changes dramatically with regard to the existence and characterization
of extremal domains. One consequence is that area-normalized eigenvalues may now remain bounded from both above and below: we prove
in \autoref{linearbound} that the scaled and normalized first eigenvalue is maximal for the degenerate rectangle, for each real $\alpha$, and for each positive $\alpha$ the eigenvalue is  bounded below (since it is positive). In \autoref{mapping} we show that if one normalizes not the area of the domain but rather its conformal mapping radius, while maintaining the perimeter scaling, then the disk is promoted to maximise the first eigenvalue. 

The above result for the first eigenvalue hints at a possible prolongation of the Szeg\H{o}--Weinberger upper bound \cite{S54,W56} for the second eigenvalue from $\alpha=0$ to $\alpha \neq 0$. The first Neumann eigenvalue is zero for all domains and so has no preferred extremal domain. The second Neumann eigenvalue is maximal for the disk by the Szeg\H{o}--Weinberger result and one hopes for this to extend to the Robin eigenvalues, at least when $|\alpha|$ is small. Indeed, for $\alpha \in [-2\pi,2\pi]$ we show in \autoref{robinszego} that the second eigenvalue is maximal for the disk among simply-connected planar domains, when the Robin parameter is scaled by perimeter and the domain is normalized by area. Hence we unify two results: Weinstock's upper bound on the first nonzero Steklov eigenvalue for domains with given perimeter, and Szeg\H{o}'s upper bound on the first nonzero Neumann eigenvalue for domains with given area. We also provide an estimate on the value of $\alpha>0$ after which the disk can no longer remain the maximal domain.

Maximality of the disk for $\lambda_{2}\left(\Omega;\alpha/L(\Omega)\right)$
implies maximality of the disk for the unscaled eigenvalue $\lambda_{2}\left(\Omega;\alpha\right)$, when $\alpha<0$, as we will show in \autoref{fixedalpha}. The point is that the unscaled eigenvalue under area normalization is equivalent to $\lambda_{2}(\Omega;\alpha/|\Omega|^{1/2})$, where the Robin parameter is scaled not by perimeter but by the square root of area, and then the (geometric) isoperimetric inequality can be applied. This corollary recovers a planar case of our earlier result that the ball maximizes the second Robin eigenvalue among domains of fixed volume \cite{FL18a}. Thus, under some circumstances, length scaling of the Robin parameter yields a stronger result than for the unscaled problem. 

Does the Rayleigh-type lower bound of Bossel \cite{B86} generalize to perimeter scaling? 
\begin{conjecture}[perimeter scaling $\Longrightarrow$ $\lambda_1$ minimal for disk]\label{ourconj}
The disk minimises 
\[
\lambda_1\left(\Omega;\frac{\alpha}{L(\Omega)}\right) 
\]
among all convex bounded planar domains of given area, when $\alpha \in \R$. 
\end{conjecture}
The restriction to convex domains is needed since D. Bucur (private communication) has pointed out counterexamples by outward boundary perturbation that drive the first eigenvalue to zero, when $\alpha>0$. 

If true, \autoref{ourconj} would imply, in the class of convex domains of given area, Bossel's result giving minimality of $\lambda_1(\Omega;\alpha)$ for the disk when $\alpha>0$. For $\alpha<0$ the conjecture is of a new and different nature, since the first eigenvalue is not even bounded below, without scaling the Robin parameter. 

In this paper we concentrate on the $2$-dimensional problem, but our proposed scaling and normalization extend naturally to the general dimension $n$ by considering quantities of the form
\[
 |\Omega|^{2/n}\lambda \! \left(\Omega;\fr{ |\Omega|^{1-2/n}}{|\partial\Omega|}\alpha \right) .
\]
The upper bound on the first eigenvalue in \autoref{linearbound} extends to higher dimensions in this manner, with an analogous proof. For maximising the second eigenvalue, we raise a higher dimensional conjecture for convex domains in \autoref{sopenprob}, where some other open problems are discussed too.

\section{\bf Notation and main results}

We consider the quantity
\[
\lambda_k\big( \Omega;\alpha/L(\Omega) \big) A(\Omega) , \qquad k=1,2,\dots ,
\]
in which each eigenvalue is multiplied by the area $A(\Omega)$, and the Robin parameter is scaled by the perimeter $L(\Omega)$. This quantity is scale invariant --- its value does not change when $\Omega$ is scaled by a positive constant factor $t$, thanks to the rescaling relation $t^2 \lambda(t\Omega;\alpha/t) = \lambda(\Omega;\alpha)$. In terms of Rayleigh quotients, the one associated to $\lambda_k(\Omega;\alpha)$ in \autoref{robinproblem} is
\[
Q[u] = 
\frac{\int_\Omega |\nabla u|^2 \, dx + \alpha \int_{\partial \Omega} u^2 \, ds}{\int_\Omega u^2 \, dx} 
\]
where $u \in H^1(\Omega)$. After multiplying by area and replacing $\alpha$ with $\alpha/L(\Omega)$, the Rayleigh quotient takes an appealing ``mean value'' form,
\[
\overline{Q}[u] = \frac{A(\Omega) \fint_\Omega |\nabla u|^2 \, dx + \alpha \fint_{\partial \Omega} u^2 \, dS}{\fint_\Omega u^2 \, dx} ,
\]
where we observe that each of the three terms is scale invariant by itself. 

The distinction between the \textbf{normalizing factor} that multiplies the eigenvalue and the \textbf{scale factor} that divides the Robin parameter is central to this paper. These two distinct factors lie behind the unification (in \autoref{szegoweinstock}) of Weinstock's bound on the first Steklov eigenvalue for given perimeter and Szeg\H{o}'s bound on the first (nontrivial) Neumann
eigenvalue for given area. 

\subsection*{The first eigenvalue} Under normalization by $A$ and scaling by $L$, the first eigenvalue is  bounded from above on general domains for all $\alpha$, being maximal in the limiting case of a degenerate rectangle. This upper bound is  
elementary, yet suggestive of the different type of results we should expect now that the Robin parameter is appropriately scaled. The theorem also has the virtue of holding for all domains, and for both positive and negative values of the parameter $\alpha$. 
\begin{thmx}[Sharp upper bound on $\lambda_1$ for all $\alpha$] \label{linearbound}
Fix $\alpha \neq 0$. If $\Omega$ is a bounded, Lipschitz planar domain then 
\[
\lambda_1\big( \Omega;\alpha/L(\Omega) \big) A(\Omega) < \alpha   
\]
with equality holding in the limit for rectangular domains that degenerate to a line segment. 
\end{thmx}
In the omitted case of vanishing $\alpha$, equality holds for all domains since $\lambda_1(\Omega;0) = 0$. 

Although sharp among all domains, the theorem is not sharp for a fixed domain $\Omega$ in the limit as $\alpha$ approaches
$\pm\infty$, in the sense that the first Robin eigenvalue for any given domain approaches a finite number (the Dirichlet eigenvalue)
as $\alpha \to +\infty$, and approaches $-\infty$ quadratically rather than linearly as $\alpha \to -\infty$, by the asymptotic formula of
Lacey \emph{et al.}\ \cite[Theorem 4.14]{H17}.

\autoref{linearbound} in the nonstrict, unscaled form $\lambda_1(\Omega;\alpha) \leq \alpha L(\Omega)/A(\Omega)$ was noted by several authors previously. The novelty here consists rather of the scaling form and the asymptotic sharpness of the strict inequality. 

As mentioned in the Introduction, in $n$ dimensions \autoref{linearbound} can be generalized in a straightforward fashion to apply to
$\lambda_1(\Omega;\alpha V^{1-2/n} /S) V^{2/n}$ where $V$ is volume and $S$ is surface area. 

\subsection*{The second eigenvalue} 

A \emph{Jordan domain} is a simply-connected, bounded planar domain $\Omega$ whose boundary is a Jordan curve.
A \emph{Jordan--Lipschitz domain} is a Jordan domain with Lipschitz boundary. 
\begin{thmx}[perimeter scaling $\Longrightarrow$ $\lambda_2$ maximal for disk] \label{robinszego}
Fix $\alpha \in [-2\pi,2\pi]$. If $\Omega$ is a Jordan--Lipschitz domain then the scale invariant quantity 
\[
\lambda_2\big( \Omega;\alpha/L(\Omega) \big) A(\Omega)
\]
is maximal for the disk. Equivalently, 
\[
\lambda_2\big( \Omega;\alpha/L(\Omega) \big) \leq \lambda_2\big( D;\alpha/L(D) \big) 
\]
where $D$ is a disk with the same area as $\Omega$. Equality holds if and only if $\Omega$ is a disk.
\end{thmx}
The endpoint value $\alpha=-2\pi$ is special, because it is where $\lambda_2\big( D;\alpha/L(D)\big)=0$; indeed, by \autoref{basic2} later
in the paper, the disk $D$ of radius $R$ and perimeter $L(D)=2\pi R$ has repeated second eigenvalue $\lambda_2(D;-1/R) = \lambda_3(D;-1/R) = 0$. The corresponding eigenfunctions are $u=x_1$ and $u=x_2$. 

The hypothesis $\alpha \geq -2\pi$ in \autoref{robinszego} could perhaps be relaxed, and in fact we believe the theorem might hold for
all negative $\alpha$, on simply-connected domains.  

Similarly, we expect the hypothesis $\alpha \leq 2\pi$ is not best possible and the theorem should hold for a  larger range of $\alpha$-values. However, in this direction we know that the theorem cannot hold at $\alpha = \infty$ due to Dirichlet eigenvalues being arbitrarily large on long
thin domains. In fact, such domains show that the theorem definitely fails for $\alpha \geq 32.7$, as explained at the end of \autoref{robinszegoproof}. 

Incidentally, the reason one may state the theorem in terms of absolute constants $\pm 2\pi$ and state the counterexample with absolute constant $32.7$ is because $\alpha$ was divided by $L(\Omega)$. Otherwise the perimeter would need to be included in all the relevant statements and constants. 

The Lipschitz assumption on the boundary in \autoref{robinszego} could be weakened somewhat, since it is used only to guarantee compactness of the imbedding $H^1 \hookrightarrow L^2$ and existence of the trace operator on the boundary, and to ensure the chord--arc condition in Case (ii) of \autoref{robinszegoproof}. 

\smallskip
A corollary with fixed negative $\alpha$ (not scaled by perimeter) follows easily from the theorem with the help of the isoperimetric inequality. Let $R(\Omega)=\sqrt{A(\Omega)/\pi}$ be the radius of the disk $D$ having the same area as $\Omega$.
\begin{corx}[no scaling $\Longrightarrow$ $\lambda_2$ maximal for disk] \label{fixedalpha}
If $\Omega$ is a Jordan--Lipschitz domain and $\alpha \in [-1/R(\Omega),0]$ then $\lambda_2(\Omega;\alpha) \leq \lambda_2(D;\alpha)$, with equality if and only if $\Omega$ is a disk.
\end{corx}
This corollary is a special case of our earlier result \cite[Theorem A]{FL18a} for arbitrary domains in all dimensions, which was proved using a Weinberger-type method. Thus for the second Robin eigenvalue on simply-connected planar domains, the Szeg\H{o} method gives a definitely stronger inequality (\autoref{robinszego}) than the Weinberger method (\autoref{fixedalpha}). On the other hand, the Weinberger method offers additional flexibility, which we exploited in \cite[Theorem A]{FL18a} to prove the result of \autoref{fixedalpha} for a larger range of $\alpha$-values, all the way down to $-3/2R(\Omega)$. Further, Weinberger's method works regardless of connectivity, whereas \autoref{robinszego} fails for certain doubly connected domains (annuli), as explained below. 

We shall now relate our results to the Neumann and Steklov spectra. To this end write $0=\mu_0 < \mu_1 \leq \mu_2 \leq \dots$ for the spectrum
of the Neumann Laplacian, and $0=\sigma_0 < \sigma_1 \leq \sigma_2 \leq \dots$ for the Steklov spectrum (corresponding to harmonic functions
with $\partial u/\partial \nu = \sigma u$ on the boundary). For an introduction to Steklov spectral geometry, we highly recommend Girouard and Polterovich's survey paper \cite{GP17}. 

The next result unifies Weinstock's upper bound on $\sigma_1$ under perimeter normalization with Szeg\H{o}'s upper bound on $\mu_1$ under area normalization. Until now these results have been regarded as different due to their different normalizing factors, although the proofs are clearly closely related \cite{GP10}. By inspecting the horizontal and vertical intercepts of $\alpha \mapsto \lambda_2(\Omega;\alpha/L)A$, we discover that the Steklov and Neumann inequalities are in fact two facets of one underlying result, \autoref{robinszego}.  
\begin{corx}[Weinstock \protect{\cite{We54}}, Szeg\H{o} \protect{\cite{S54}}] \label{szegoweinstock}
For $\Omega$ a Jordan--Lipschitz domain, the scale invariant quantities
\[
\sigma_1(\Omega)L(\Omega) \quad \text{and} \quad \mu_1(\Omega) A(\Omega) 
\]
are maximal for the disk, and only for the disk. 
\end{corx}

The Weinstock inequality on $\sigma_1(\Omega)L(\Omega)$ fails for certain annuli \cite[Example 5.14]{H17}. Hence the above corollary
and \autoref{robinszego} both fail for general domains that are not simply connected. On the other hand, by weakening the normalization to area and considering $\sigma_1(\Omega)\sqrt{A(\Omega)}$, Brock did obtain a result valid for all domains, and which extends to all dimensions~ \cite{B01}.
The Szeg\H{o} inequality on $\mu_1(\Omega) A(\Omega)$ likewise holds for all domains and extends to all dimensions, as was shown by
Weinberger \cite{W56}. These Brock and Weinberger inequalities are unified by our recent work on the Robin spectrum under volume normalization with no scaling of the Robin parameter \cite[Corollary~B]{FL18a}. 

\subsection*{Other normalizations}

\smallskip
If instead of normalizing the Robin eigenvalue with area we normalize with the square of the conformal mapping radius, then for positive $\alpha$ a geometrically sharp result can be obtained for the first eigenvalue. The Robin parameter continues to be scaled by perimeter in what follows.  
\begin{thmx}[conformal radius normalization $\Longrightarrow$ $\lambda_1$ maximal for the disk] \label{mapping}
Suppose $F : \D \to \Omega$ is a conformal map of the unit disk onto a Jordan--Lipschitz domain $\Omega$. If $\alpha > 0$ then the scale invariant quantity
\[
\lambda_1\big( \Omega;\alpha/L(\Omega) \big) |F^\prime(0)|^2 
\]
is maximal if and only if $F$ is linear and $\Omega$ is a disk. 
\end{thmx}
By letting $\alpha \to \infty$, one recovers the result of P\'{o}lya and Szeg\H{o} \cite[{\S}5.8]{PS51} that the first Dirichlet eigenvalue normalized by conformal mapping radius, $\lambda_1^\text{Dir} |F^\prime(0)|^2$, is maximal for the disk.

\section{\bf Open problems and conjectures\label{sopenprob}}

A stronger result than \autoref{szegoweinstock} is known to hold, namely, that the normalized harmonic means 
\[
\frac{L}{(\sigma_1^{-1} + \sigma_2^{-1})/2} \qquad \text{and} \qquad \frac{A}{(\mu_1^{-1} + \mu_2^{-1})/2}
\]
of the first two Steklov and Neumann eigenvalues are maximal for the disk, among simply-connected domains; see \cite[p.\ 634]{W56}.
A natural question is whether \autoref{robinszego} can be strengthened in a similar way to handle the harmonic mean of the Robin
eigenvalues $\lambda_2$ and $\lambda_3$.

Another open problem is to generalize \autoref{robinszego} to higher dimensions, where convexity might provide a reasonable substitute for simply connectedness. Given a domain $\Omega$ in higher dimensions, write $V$ for its volume and $S$ for its surface area. Let $\B$ be the unit ball. 
\begin{conjecture}[perimeter-volume scaling $\Longrightarrow$ $\lambda_2$ maximal for ball]\label{bucurconj}
The ball maximises the scale invariant quantity
\[
\lambda_2(\Omega;\alpha V^{1-2/n}/S) V^{2/n}
\]
among all convex bounded domains in $\Rn$, when $\alpha \in [-S(\B)/V(\B)^{1-2/n},0]$. 

Consequently $\lambda_2(\Omega;\alpha) \leq \lambda_2(B;\alpha)$ for all $\alpha \in [-1/R,0]$, where $B=B(R)$ is a ball having the same volume as $\Omega$. 
\end{conjecture}
Taking $n=2$ reduces the conjecture back to $\lambda_2(\Omega;\alpha/L)A$, as in \autoref{robinszego}.

Maximality of the ball among convex domains for the normalized Steklov eigenvalue $\sigma_1 S/V^{1-2/n}$ would follow from \autoref{bucurconj}, by arguing as in the plane for \autoref{szegoweinstock}. In fact, this maximality of the Steklov eigenvalue at the ball, among convex domains, has been proved directly already by Bucur \emph{et al.}\ \cite{BFNT17}, and one would like to extend their method to the Robin eigenvalue in order to prove \autoref{bucurconj}.

Does \autoref{mapping} hold also for the second Robin eigenvalue? It does in the limit $\alpha \to \infty$, because Ashbaugh and Benguria \cite[\S4]{AB93} proved for the second Dirichlet eigenvalue that $\lambda_2^\text{Dir} |F^\prime(0)|^2$ is maximal for the disk. Curiously, this result was not proved by employing conformal mapping to create trial functions for the second eigenvalue. Instead, they combined their sharp PPW inequality on the ratio of the first two eigenvalues with P\'{o}lya and Szeg\H{o}'s bound on the first eigenvalue, using the decomposition
\[
\lambda_2^\text{Dir} |F^\prime(0)|^2 = \frac{\lambda_2^\text{Dir}}{\lambda_1^\text{Dir}} \left( \lambda_1^\text{Dir} |F^\prime(0)|^2 \right)
\]
where each factor on the right side is maximal for the disk. In view of this Dirichlet result it seems natural to conjecture that the second Robin eigenvalue is maximised at the disk.
\begin{conjecture}[conformal radius normalization $\Longrightarrow$ $\lambda_2$ maximal for the disk] \label{mapping2}
Suppose $F : \D \to \Omega$ is a conformal map of the unit disk onto a Jordan--Lipschitz domain $\Omega$. If $\alpha > 0$ then the scale invariant quantity
\[
\lambda_2\big( \Omega;\alpha/L(\Omega) \big) |F^\prime(0)|^2 
\]
is maximal when $F$ is linear and $\Omega$ is a disk. 
\end{conjecture}
We already discussed the limit $\alpha \to \infty$. At the other extreme, when $\alpha=0$ the conjecture says $\mu_1 |F^\prime(0)|^2$ is maximal when $F$ is linear and $\Omega$ is the disk, where $\mu_1$ is the first positive eigenvalue of the Neumann Laplacian. This claim is certainly true, as it follows from Szeg\H{o}'s theorem \cite{S54} maximising $\mu_1 A$ for the disk, noting that the ratio $|F^\prime(0)|^2/A=|F^\prime(0)|^2/\int_\D |F^\prime(z)|^2 \, |dz|^2$ is maximal when $F^\prime$ is constant, that is, when $F$ is linear and $\Omega$ is a disk.

\autoref{mapping} could perhaps be generalized to cone metrics on the disk and other geometric situations considered in the Dirichlet case by Laugesen and Morpurgo \cite{LM98}.

\subsection*{Eigenvalue sums} The methods of this paper do not seem to extend to eigenvalue sums of the form $\lambda_1 + \dots + \lambda_m$, because composition with a conformal map does not preserve $L^2$-orthogonality of trial functions, while pre-composition with a M\"{o}bius transformation of the disk can help only to the extent of a few degrees of freedom. 

Composition with a linear transformation, on the other hand, does preserve $L^2$-orthogonality. That observation has generated a number of sharp upper bounds on sums of Robin and magnetic Robin eigenvalues for domains that are linear images of rotationally symmetric domains, in work by Laugesen \emph{at al.}\ \cite[Theorem 3.2]{LLR12} and Laugesen and Siudeja \cite[Theorem 3.3]{LS11a}, \cite[Theorem 3]{LS11b}, with generalizations to starlike domains also \cite[Theorem 3.5]{LS14}. The Robin parameter in these results is scaled by various geometric factors of the domain such as its moment of inertia \cite[Lemma 5.3]{LS11a}, and thus the scaling is more complicated than the perimeter factor used in this paper. 

The methods of this paper also do not appear to extend to reciprocal sums of the form $1/\lambda_1 + \dots + 1/\lambda_m$ or to spectral zeta functions, because the numerator $\int_\Omega |\nabla u|^2 \, dx + \alpha \int_{\partial \Omega} u^2 \, ds$ of the Robin Rayleigh quotient is not conformally invariant. 

\subsection*{Lower bounds --- literature and discussion} To complete the context for the current paper's upper bounds on eigenvalues, we mention the Faber--Krahn type lower bound on the first eigenvalue, $\lambda_1(\Omega;\alpha)A(\Omega) \geq \lambda_1(D;\alpha)A(D)$, proved for $\alpha>0$ by Bossel \cite{B86} and extended to the $n$-dimensional case by Daners \cite{D06}. An alternative approach via the calculus of variations was found more recently by Bucur and Giacomini \cite{BG10,BG15a}, with a quantitative version by Bucur \emph{et al.}\ \cite{BFNT18}. Among the family of rectangles of given area the square is the minimizer \cite[Theorem 4.1]{FK18}, with a generalization to higher dimensions through appealing convexity arguments by Keady and Wiwatanapataphee \cite{KW18}. 

For the reverse inequality when $\alpha<0$, which is known as the Bareket conjecture, a great deal is now known for domains near the disk by Ferone \emph{et al.}\ \cite{FNT15}, and for general domains when $|\alpha|$ is small by Freitas and Krej\v{c}i\v{r}\'{\i}k \cite{FK15}, while annular counterexamples have been discovered for large $|\alpha|$. References and a fuller discussion are provided in our earlier paper \cite[{\S}1]{FL18a}.

For a lower bound on the second eigenvalue, Kennedy \cite{K09} observed that Krahn's two-disk argument for the Dirichlet Laplacian carries across to the Robin case as a corollary of Bossel's inequality for the first eigenvalue. For more on spectral shape optimization we recommend the survey volume edited by Henrot \cite{H17}. 

\section{\bf Proof of \autoref{linearbound}}
\label{linearboundproof}
Substituting the constant trial function $u(x) \equiv 1$ into the Rayleigh quotient gives the upper bound
\[
   \lambda_1(\Omega;\alpha/L) A \leq \frac{0 + (\alpha/L)\int_{\partial \Omega} 1^2 \, ds}{\int_\Omega 1^2 \, dx} \, A = \alpha.
\]
We show this inequality must be strict. If equality held, then the constant trial function $u$ would be a first eigenfunction, and so $\lambda_1(\Omega;\alpha/L) u = - \Delta u = 0$, which means $\lambda_1(\Omega;\alpha/L)=0$. From equality holding we would deduce $\alpha=0$, contradicting a hypothesis in the theorem. Hence equality cannot hold and the inequality is strict. 

To show equality is attained asymptotically for rectangles degenerating to a line segment, consider the family of rectangles $\Omega_t$ having side lengths $t$ and $1/t$, area $A(t)=1$ and perimeter $L(t)=2(t+t^{-1})$, where $t \geq 1$. By separation of variables and using a known lower bound on the first eigenvalue of an interval \cite[Appendix A.1]{FK18}, one gets for fixed $\alpha>0$ that 
\[
 \lambda_1\big(\Omega_t;\alpha/L(t)\big)A(t) \geq \alpha - O_\alpha(t^{-2}) \qquad \text{as $t \to \infty$.}
\]
Hence
\begin{equation} \label{asympineq}
\lambda_1\big(\Omega_t;\alpha/L(t)\big)A(t) \to \alpha \qquad \text{as $t \to \infty$,}
\end{equation}
and so equality is attained asymptotically in the theorem. 

The argument is similar when $\alpha<0$, by using hyperbolic trigonometric instead of trigonometric functions for the separated eigenfunctions. 
\qed

\section{\bf The Robin spectrum on the disk}
\label{preliminaries}

The proof of \autoref{robinszego} will require some properties of the Robin eigenvalues and eigenfunctions on the unit disk $\D$. Separating variables in the Robin eigenvalue problem \eqref{robinproblem} with   
\[
u(r,\theta)=g(r)T(\theta) 
\] 
implies that the angular part satisfies $T^{\prime \prime}(\theta) + \kappa^2 T(\theta) = 0$ 
where $\kappa \geq 0$ is an integer. When $\kappa=0$ (giving a constant function $T$) the 
eigenfunctions on the disk are purely radial. For positive values of $\kappa$ the  angular function $T(\theta)$ equals $\cos \kappa \theta$ or $\sin \kappa \theta$, and the eigenvalues have multiplicity $2$.

\begin{figure}
\includegraphics[scale=0.6]{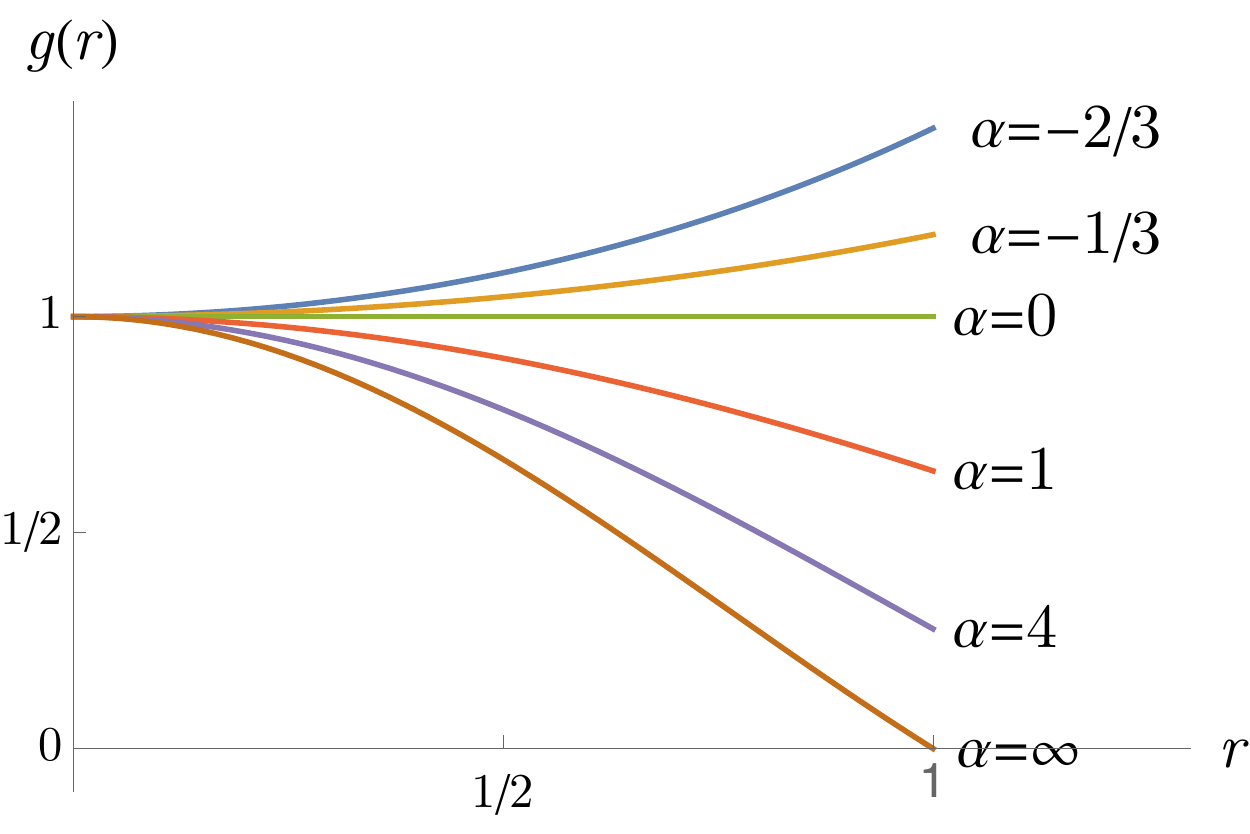}
\caption{\label{Robin_g_first}Plot of the first Robin eigenfunction $g(r)$ of the unit disk, for various values of $\alpha$, normalized with $g(0)=1$. When $\alpha=0$ one sees $g(r)$ is the constant Neumann eigenfunction with eigenvalue $0$, and when $\alpha=\infty$ it is the Dirichlet eigenfunction $J_0(j_{0,1}r)$ with eigenvalue $j_{0,1}^2$. Between these extremes, $g(r)=J_0(\sqrt{\lambda_1}\,r)$ where $\lambda_1=\lambda_1(\D;\alpha)>0$ is the eigenvalue. }
\end{figure}

\begin{figure}
\includegraphics[scale=0.6]{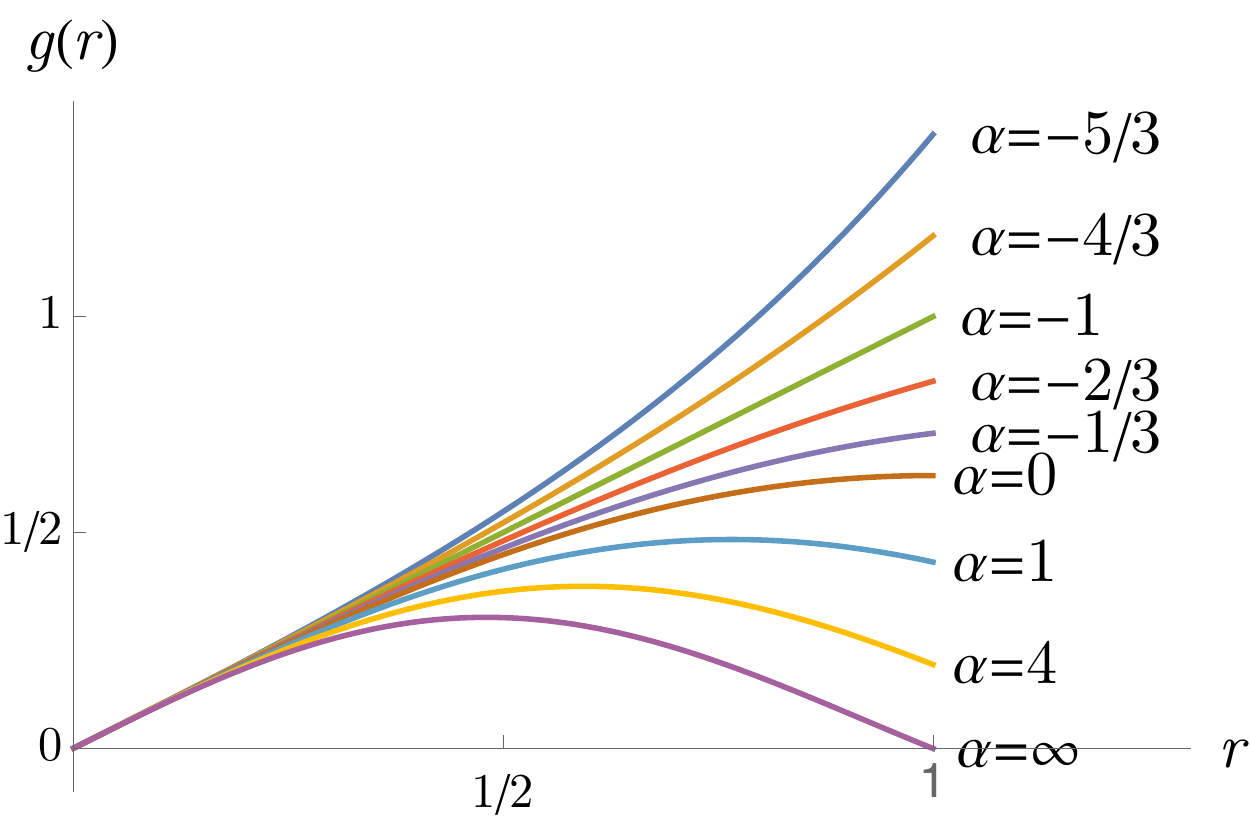}
\caption{\label{Robin_g_second}Plot of the radial part $g(r)$ of the second Robin eigenfunction of the unit disk, for various values of $\alpha$, normalized with $g^\prime(0)=1$. When $\alpha=-1$ it is the straight line $g(r)=r$ and $\lambda_2(\D;-1)=0$. When $\alpha>-1$ one has  $g(r)=(\text{const.})J_1(\sqrt{\lambda_2}\,r)$ where $\lambda_2=\lambda_2(\D;\alpha)>0$ is the eigenvalue. The eigenfunctions are $g(r) \cos \theta$ and $g(r) \sin \theta$.}
\end{figure}

The radial part $g$ satisfies the Bessel-type equation
 \[
 g''(r) + \fr{1}{r}g'(r) + \left( \lambda - \frac{\kappa^2}{r^2} \right) g(r) = 0
\]
due to the eigenfunction equation $-\Delta u = \lambda u$, while the boundary condition 
\[
\frac{\partial u}{\partial \nu} + \alpha u = 0
\]
at $r=1$ implies 
\[
g^\prime(1) + \alpha g(1) = 0 .
\] 
The key facts about the first and second eigenvalues and eigenfunctions are summarized in the next propositions and in \autoref{Robin_g_first} and \autoref{Robin_g_second}, which are taken from  \cite[Section 5]{FL18a}, where the ball was handled in all dimensions. The spectral curves for the disk are illustrated in \cite[Figure~3]{FL18a}.

For simplicity, since the domain is fixed in this section, we do not rescale $\alpha$ by the perimeter $2\pi$ of the disk. Thus the range $\alpha \in [-2\pi,2\pi]$ in \autoref{robinszego} corresponds here to $\alpha \in [-1,1]$.
\begin{proposition}[First Robin eigenfunction of the disk]\label{basic1} The first eigenvalue of $\D$ is simple, and changes sign at $\alpha=0$ according to 
\[
\lambda_1(\D;\alpha)
\begin{cases}
< 0 & \text{when $\alpha < 0$,} \\
= 0 & \text{when $\alpha = 0$,} \\
> 0 & \text{when $\alpha > 0$.}
\end{cases}
\]
The first eigenfunction is radial ($\kappa=0$), with $g(0) > 0$ and $g^\prime(0)=0$. If $\alpha< 0$ then $g^\prime(r) > 0$; if $\alpha=0$ then $g^\prime(r) = 0$; and if $\alpha > 0$ then $g^\prime(r) < 0$, when $r \in (0,1)$.
\end{proposition}
\begin{proposition}[Second Robin eigenfunctions of the disk]\label{basic2} The eigenfunctions for the double eigenvalue $\lambda_2(\D;\alpha)=\lambda_3(\D;\alpha)$ have simple angular dependence ($\kappa=1$), meaning they take the form 
\[
g(r)\cos \theta \qquad \text{and} \qquad g(r) \sin \theta .
\]
The radial part has $g(0)=0,g^\prime(0)>0,g(r)>0$ for $r \in (0,1)$, and $g(1)>0$. When $\alpha \leq 0$ one finds $g(r)$ is strictly increasing, with $g^\prime(r)>0$. When $\alpha > 0$, the derivative $g^\prime$ is positive on some interval $(0,r_\alpha)$ and negative on $(r_\alpha,1)$, for some number $r_\alpha \in (0,1)$. 

The eigenvalue changes sign at $\alpha=-1$, with 
\[
\lambda_2(\D;\alpha)=\lambda_3(\D;\alpha)
\begin{cases}
< 0 & \text{when $\alpha < -1$,} \\
= 0 & \text{when $\alpha = -1$,} \\
> 0 & \text{when $\alpha > -1$.}
\end{cases}
\]
\end{proposition}
A couple of the assertions in \autoref{basic2} when $\alpha>0$ are not included in \cite[Section 5]{FL18a}, and so we justify them here. The radial part of the second Robin eigenfunction is $g(r)=(\text{const.})J_1(\sqrt{\lambda_2(\D;\alpha)}\,r)$. As $\alpha$ increases from $0$ to $\infty$ the eigenvalue increases from the second Neumann eigenvalue to the second Dirichlet eigenvalue of the unit disk, and so $j_{1,1}^\prime < \sqrt{\lambda_2(\D;\alpha)} < j_{1,1}$. (Numerically, $j_{1,1}^\prime \simeq 1.84$ and $j_{1,1} \sim 3.83$.) The Bessel function $J_1$ vanishes at $0$ and at $j_{1,1}$, and has positive derivative on $(0,j_{1,1}^\prime)$ and negative derivative on $(j_{1,1}^\prime,j_{1,1})$. Hence $g(1)>0$, and $g^\prime$ is positive on the interval $(0,r_\alpha)$ and negative on $(r_\alpha,1)$, where the number $r_\alpha = j_{1,1}^\prime/\sqrt{\lambda_2(\D;\alpha)}$ lies between $0$ and $1$.

\section{\bf Center of mass argument}
\label{centerofmasssection}

In this section, $\Omega$ is a simply-connected planar domain and $g(r)$ is a continuous function for $0 \leq r \leq 1$ with $0=g(0)<g(1)$. Define continuous functions
\begin{equation*} 
u_2 = g(r) \cos \theta , \qquad u_3 = g(r) \sin \theta ,
\end{equation*}
on the unit disk $\D$. The following center of mass result will be used in proving \autoref{robinszego}. 
\begin{lemma}[Center of mass]\label{centerofmass}
If $v_1$ is an integrable real-valued function on $\Omega$ with $\int_\Omega v_1 \, dx > 0$, then a conformal map $f : \D \to \Omega$ can be chosen such that the functions $v_2 = u_2 \circ f^{-1}$ and $v_3 = u_3 \circ f^{-1}$ are orthogonal to $v_1$:
\[
\int_\Omega v_2 v_1 \, dx = 0 \qquad \text{and} \qquad \int_\Omega v_3  v_1\, dx = 0 .
\]
\end{lemma}
Szeg\H{o} \cite[Section 2.5]{S54} treated the case $v_1 \equiv 1$ by an approximate identity argument and elementary index theory. Hersch \cite{H70} reformulated the argument more geometrically, avoiding Szeg\H{o}'s use of approximate identities. For the sake of completeness, we include a version of Hersch's proof below. 

\begin{proof}[Proof of \protect{\autoref{centerofmass}}] Fix a conformal map $F : \D \to \Omega$, and let $H(z) = g(r)e^{i\theta}$ where $z=re^{i\theta} \in \overline{\D}$. Note $H$ is continuous on the closed disk, including at the origin since $g(0)=0$. Define a complex-valued function (vector field) on the disk by
\[
V(\zeta) = \int_\Omega H \big( M_\zeta\big(F^{-1}(x)\big) \big) v_1(x) \, dx , \qquad \zeta \in \D , 
\]
where 
\[
M_\zeta(z) = \frac{z+\zeta}{1+z\overline{\zeta}} , \qquad z \in \D , 
\]
is a M\"{o}bius map of the unit disk $\D$ to itself. 

Notice $M_\zeta(z)$ remains continuous as a function of $(\zeta,z) \in \overline{\D} \times \D$ (where now we allow $|\zeta|=1$), taking values in $\overline{\D}$. Thus the vector field $V(\zeta)$ is well defined for $\zeta \in \overline{\D}$, and is continuous at each point by a simple application of dominated convergence, using continuity and boundedness of $H$. The boundary behavior is easily determined: when $\zeta=e^{i\phi}$ one has $M_\zeta(z)=e^{i\phi}$ for all $z \in \D$, and so 
\[
V(e^{i\phi}) = g(1)e^{i\phi} \int_\Omega v_1 \, dx , \qquad \phi \in [0,2\pi] .
\]
Thus the continuous vector field $V$ points radially outward on the unit circle,  because $g(1) \int_\Omega v_1 \, dx > 0$ by construction. 

Index theory, or the Brouwer fixed point theorem, implies that $V$ vanishes somewhere in the interior of the disk. That is, $V(\zeta)=0$ for some $\zeta \in \D$, which means $H \circ f^{-1}$ is orthogonal to $v_1$, where $f = F \circ M_\zeta^{-1}$. Because $H=u_2+iu_3$ by definition, we conclude $u_2 \circ f^{-1}$ and $u_3 \circ f^{-1}$ are orthogonal to $v_1$. 
\end{proof}

\section{\bf Proof of \autoref{robinszego}}
\label{robinszegoproof}
After rescaling, we may suppose $\Omega$ has area $\pi$, so that $D$ is the unit disk $\D$. Our goal is to show
\[
\lambda_2\big( \Omega;\alpha/L(\Omega) \big) \leq \lambda_2\big( \D;\alpha/2\pi \big) , \qquad \alpha \in [-2\pi,2\pi] .
\]
Note on the right side that $\lambda_2\big(\D;\alpha/2\pi \big) \geq 0$ by \autoref{basic2}, since $\alpha/2\pi \geq -1$. Thus if $\lambda_2\big(\Omega;\alpha/L(\Omega)\big) < 0$ then there is nothing to prove. 

\subsection*{Case (i). Second eigenvalue greater than zero} Assume $\lambda_2\big(\Omega;\alpha/L(\Omega)\big) > 0$. Let $u_2$ and $u_3$ be the second Robin eigenfunctions of the unit disk with Robin parameter $\alpha/2\pi$, which from \autoref{basic2} have the form
\[
u_2 = g(r) \cos \theta \qquad \qquad  \text{and} \qquad u_3 = g(r) \sin \theta ,
\]
where $g$ is smooth with $0=g(0)<g(1)$. Take a conformal map $f$ from $\D$ onto $\Omega$, and define  
\[
v_2 = u_2 \circ f^{-1} \qquad \text{and} \qquad v_3 = u_3 \circ f^{-1} .
\]
These functions belong to $H^1(\Omega)$ because they are bounded and smooth with 
\begin{equation} \label{confinv}
\int_\Omega |\nabla v_k|^2 \, dx = \int_\D |\nabla u_k|^2 \, dx < \infty, \qquad k=2,3,
\end{equation}
by conformal invariance of the Dirichlet integral. Note $v_2$ and $v_3$ extend continuously to $\partial \Omega$ since $f^{-1}$ extends continuously (using that $\partial \Omega$ is a Jordan curve). 

The conformal map can be chosen by \autoref{centerofmass} to ensure the orthogonality relations
\[
\int_\Omega v_2 v_1 \, dx = 0 \qquad \text{and} \qquad \int_\Omega v_3  v_1\, dx = 0 ,
\]
where $v_1$ is the first Robin eigenfunction on $\Omega$ for Robin parameter $\alpha/L(\Omega)$; note here that $v_1$ does not change sign and so we may assume its integral is positive. Thus $v_2$ and $v_3$ are valid trial functions for $\lambda_2(\Omega;\alpha/L(\Omega))$.  Taking $v_2$ as a trial function in the Rayleigh principle for the second eigenvalue shows that
\[
\lambda_2\big(\Omega;\alpha/L(\Omega)\big) \int_\Omega v_2^2 \, dx 
\leq \int_\Omega |\nabla v_2|^2 \, dx + \frac{\alpha}{L(\Omega)} \int_{\partial \Omega} v_2^2 \, ds . 
\]
This formula pulls back under the conformal map to 
\[
\lambda_2\big(\Omega;\alpha/L(\Omega)\big) \int_\D u_2^2 |f^\prime|^2 \, dx 
\leq \int_\D |\nabla u_2|^2 \, dx + \frac{\alpha}{L(\Omega)} \int_{\partial \Omega} v_2^2 \, ds ,
\]
due to the conformal invariance in \eqref{confinv}. Substituting the definition $u_2 = g(r) \cos \theta$ gives
\[
\begin{split}
& \lambda_2\big(\Omega;\alpha/L(\Omega)\big) \int_\D g(r)^2 (\cos \theta)^2 |f^\prime|^2 \, dx \\
\leq & \int_\D \left( g^\prime(r)^2 \cos^2 \theta + r^{-2} g(r)^2 \sin^2 \theta \right) \, dx + \frac{\alpha}{L(\Omega)} \int_{\partial \Omega} v_2^2 \, ds .
\end{split}
\]

An analogous formula holds for $u_3$, with the roles of $\cos$ and $\sin$ interchanged. Adding that formula to the preceding one and using that $\cos^2 \theta + \sin^2 \theta = 1$ and hence  $v_2^2+v_3^2=g(1)^2$ on $\partial \Omega$, we deduce
\begin{equation} \label{mainformula}
\lambda_2\big(\Omega;\alpha/L(\Omega)\big) \int_\D g(r)^2 |f^\prime|^2 \, dx 
\leq \int_\D \left( g^\prime(r)^2 + r^{-2} g(r)^2 \right) \, dx + \alpha g(1)^2 .
\end{equation}
Equality holds if $\Omega$ is the unit disk and $f$ is the identity map, since $u_2$ and $u_3$ are  the second eigenfunctions of the disk, which means 
\begin{equation} \label{mainformulaequality}
\lambda_2(\D;\alpha/2\pi) \int_\D g(r)^2 \, dx 
= \int_\D \left( g^\prime(r)^2 + r^{-2} g(r)^2 \right) \, dx + \alpha g(1)^2 .
\end{equation}

Suppose for the remainder of Case (i) that $\Omega$ is not a disk. We will show  
\begin{equation} \label{areabound}
\int_\D g(r)^2 \, dx < \int_\D g(r)^2 |f^\prime|^2 \, dx .
\end{equation}
This estimate \eqref{areabound} can then be substituted into the left side of \eqref{mainformula}, relying here on the positivity of $\lambda_2$. Combining the resulting inequality with equality \eqref{mainformulaequality} for the disk, once concludes that $\lambda_2\big(\Omega;\alpha/L(\Omega)\big) < \lambda_2\big(\D;\alpha/2\pi\big)$, as wanted for the theorem.

Szeg\H{o} proved inequality \eqref{areabound} under the assumption that $g$ is increasing, which in our Robin situation holds when $\alpha \leq 0$. We will extend his method to handle $\alpha \leq 2\pi$. To start with, 
\begin{align*}
\int_\D g(r)^2 |f^\prime|^2 \, dx - \int_\D g(r)^2 \, dx
& = - \int_0^1 g(r)^2 \frac{d\ }{dr} \left( \pi r^2 - \int_{\D(r)} |f^\prime|^2 \, dx \right) dr \\
& = \int_0^1 2g(r)g^\prime(r) \left( \pi r^2 - \int_{\D(r)} |f^\prime|^2 \, dx \right) dr 
\end{align*}
by integration by parts, noting that the boundary terms vanish because 
\begin{equation} \label{mvR}
\pi = A(\Omega) = \int_\D |f^\prime|^2 \, dx .
\end{equation}
Hence 
\begin{equation} \label{areabound2}
\int_\D g(r)^2 |f^\prime|^2 \, dx - \int_\D g(r)^2 \, dx
= \int_0^1 2g(r)g^\prime(r) \pi r^2 \left( 1 - M(r) \right) dr 
\end{equation}
where 
\[
M(r) = \frac{1}{\pi r^2 } \int_{\D(r)} |f^\prime|^2 \, dx 
\]
is the mean value function. The mean value is increasing due to subharmonicity of $|f^\prime|^2$. More directly, one may write $f(z)=\sum_{n=0}^\infty a_n z^n$ as a power series and substitute into $M(r)$ to obtain 
\begin{equation} \label{mvp}
M(r) = \frac{1}{\pi r^2}\int_{\D(r)} |f^\prime(\rho e^{i\theta})|^2 \, \rho \, d\rho d\theta = \sum_{n=1}^\infty n|a_n|^2 r^{2(n-1)} ,
\end{equation}
which plainly increases as a function of $r$. Further, since $\Omega$ is not a disk we have $f(z) \not\equiv a_0+a_1 z$ and so $a_n \neq 0$ for some $n \geq 2$, which implies by \eqref{mvp} that $M(r)$ is strictly increasing as a function of $r$. 

The area normalization \eqref{mvR} gives $M(1)=1$, and so $M(r) < 1$ for $r \in (0,1)$.

If $\alpha \leq 0$ then $g$ and $g^\prime$ are both positive on $(0,1)$ by \autoref{basic2}, and so  inequality \eqref{areabound} follows from \eqref{areabound2}.

Next assume $0 < \alpha \leq 2\pi$. Define
\[
G(r) = \int_0^r 2g(\rho)g^\prime(\rho) \pi \rho^2 \, d\rho , \qquad r \in [0,1] . 
\]
Formula \eqref{areabound2} becomes
\begin{align*}
\int_\D g(r)^2 |f^\prime|^2 \, dx - \int_\D g(r)^2 \, dx
& = \int_0^1 G^\prime(r) \left( 1 - M(r) \right) dr \\
& = \int_0^1 G(r) M^\prime(r) \, dr 
\end{align*}
after integrating by parts, since $G(0)=0$ and $M(1)=1$. We want this last integral to be positive, so that \eqref{areabound} holds. Because $M^\prime(r) > 0$, it suffices to show $G(r) > 0$ for $r \in (0,1)$. Recall from \autoref{basic2} that when $\alpha > 0$, the function $g$ is positive on $(0,1)$ while $g^\prime$ is positive on some interval $(0,r_\alpha)$ and negative on $(r_\alpha,1)$. Thus $G^\prime$ is positive on $(0,r_\alpha)$ and negative on $(r_\alpha,1)$, and so to show $G(r)$ is positive for $r \in (0,1)$, we need only show $G(1) \geq 0$. That is, we want 
\[
\int_0^1 2g(r)g^\prime(r) \pi r^2 \, dr \geq 0.
\]

By \autoref{preliminaries} one has $g(r)=J_1(\sqrt{\lambda_2}r)$, where $\lambda_2=\lambda_2\big(\D;\alpha/2\pi \big)>0$. Applying this formula for $g$ and making a change of variable, we reduce to showing 
\[
\int_0^{\sqrt{\lambda_2}} 2 J_1(r)J_1^\prime(r) \, r^2 \, dr \geq 0 .
\]
The antiderivative for the left side is $J_0(r) J_2(r) r^2$, as one can check using standard Bessel formulas \cite[Eq.~(10.6.1) and (10.6.2)]{DLMF}. Hence the inequality to be proved is 
\[
J_0(\sqrt{\lambda_2}) J_2(\sqrt{\lambda_2}) \lambda_2 \geq 0 .
\] 
Note the second Robin eigenvalue $\lambda_2$ of the disk is less than the second Dirichlet eigenvalue $j_{1,1}^2$ of the disk, which in turn is less than $j_{2,1}^2$. Therefore $J_2(\sqrt{\lambda_2}) > 0$, and so the last displayed inequality holds if and only if $J_0(\sqrt{\lambda_2}) \geq 0$. Thus we want to show $\sqrt{\lambda_2} \leq j_{0,1}$, or $\lambda_2\big(\D;\alpha/2\pi \big) \leq j_{0,1}^2$.

The Robin eigenvalue increases with $\alpha$, and since $\alpha \leq 2\pi$ by hypothesis, it suffices to take $\alpha=2\pi$ and show $\lambda_2(\D;1) \leq j_{0,1}^2$. For this, observe that $u=J_1(j_{0,1}r)\cos \theta$ is a nonradial eigenfunction of the Laplacian on the unit disk, with eigenvalue $j_{0,1}^2$. We confirm that $u$ satisfies the Robin boundary condition with $\alpha=1$, namely $\partial u/\partial \nu + u = 0$ at $r=1$, by computing  
\[
j_{0,1}J_1^\prime(j_{0,1}) + J_1(j_{0,1}) = - j_{0,1}J_0^{\prime \prime}(j_{0,1}) - J_0^\prime(j_{0,1}) = 0 ,
\]
where we used the relation $J_1=-J_0^\prime$ and the Bessel equation $r^2J_0^{\prime \prime}(r)+rJ_0^\prime(r)+r^2J_0(r)=0$. Case (i) of the proof is finished. 

\subsection*{Case (ii). Second eigenvalue equal to zero} We must still handle the situation where the second eigenvalue on $\Omega$ equals zero, that is, $\lambda_2\big(\Omega;\alpha/L(\Omega)\big) = 0$. The conclusion of the theorem is then immediate, since  
\[
\lambda_2\big(\Omega;\alpha/L(\Omega)\big) = 0 \leq \lambda_2\big(\D;\alpha/2\pi \big)
\]
by \autoref{basic2}, since $\alpha/2\pi \geq -1$.  

We will show that if equality holds, then $\Omega$ is a disk. This part of the argument follows Weinstock's equality case \cite[{\S}3]{We54}. He assumed the boundary of $\Omega$ to be analytic whereas we assume only Lipschitz smoothness. We invoke subtle results from complex analysis to ensure that the harmonic function $\log |f^\prime|$ equals the Poisson integral of its boundary values. The need for such care in the Lipschitz case may not have been recognized in earlier treatments  \cite[Theorem 1.3]{GP10}.

Suppose equality holds above, meaning $\lambda_2\big(\D;\alpha/2\pi \big) = 0$. Then $\alpha=-2\pi$ and the eigenfunctions $u_2$ and $u_3$ for the disk are the coordinate functions $x_1$ and $x_2$, by the case ``$\alpha=-1$'' in \autoref{basic2}, with $g(r)=r$.

Equality holds in \eqref{mainformula}, because both sides of the inequality equal $0$. By the Rayleigh principle, the trial functions $v_2$ and $v_3$ used to derive \eqref{mainformula} must therefore be eigenfunctions on $\Omega$ with eigenvalue $0$. Thus $v_2$ and $v_3$ satisfy the (weak form of) the Robin boundary condition, which we proceed to investigate. 

Since $\partial \Omega$ is a rectifiable Jordan curve, the derivative $f^\prime$ of the conformal map belongs to the analytic Hardy space     
and the boundary values $f^\prime(e^{i\theta})$ provide the Jacobian factor for arclength (\cite[Theorem 3.12]{D00} and remarks following it). That is, 
\[
ds = |f^\prime(e^{i\theta})| \, d\theta
\]
where $ds$ denotes the arclength element on $\partial \Omega$. We will show $|f^\prime(e^{i\theta})|$ is constant a.e.

The weak formulation of the eigenfunction equation for $v_2$ on $\Omega$, with $\alpha=-2\pi$ and eigenvalue $0$ as above, says
\[
\int_\Omega \nabla v_2 \cdot \nabla \psi \, dx - \frac{2\pi}{L(\Omega)} \int_{\partial \Omega} v_2 \psi \, ds = 0, \qquad \psi \in H^1(\Omega) .
\]
Pulling back to $\D$, we deduce by conformal invariance that
\[
\int_\D \nabla u_2 \cdot \nabla \phi \, dx = \frac{2\pi}{L(\Omega)} \int_{\partial \D} u_2 \phi |f^\prime| \, d\theta , \quad \phi \in C^\infty(\overline{\D}) ,
\]
where we note that $\psi=\phi \circ f^{-1}$ belongs to $H^1(\Omega)$. Recall $u_2=x_1=r\cos \theta$. By applying Green's theorem on the left side of the last equation, we find 
\[
\int_0^{2\pi} (\cos \theta) \phi(e^{i\theta}) \, d\theta = \frac{2\pi}{L(\Omega)} \int_0^{2\pi} (\cos \theta) \phi(e^{i\theta}) |f^\prime(e^{i\theta})| \, d\theta , \quad \phi \in C^\infty(\overline{\D}) . 
\]
Since $\phi$ is arbitrary, it follows that 
\[
\cos \theta = \frac{2\pi}{L(\Omega)} (\cos \theta) |f^\prime(e^{i\theta})|
\]
for almost every $\theta$, which means $|f^\prime(e^{i\theta})| = L(\Omega)/2\pi$ a.e. Thus $|f^\prime|$ is constant a.e.\ on the unit circle. 

We will show $|f^\prime|$ is constant on the unit disk. We start by proving the Jordan curve $J=\partial \Omega$ has the chord--arc property, meaning 
\[
\text{length}\, \big( J(x,y) \big) \leq C|x-y| , \qquad x,y \in J ,
\]
for some constant $C$, where $J(x,y)$ is the shorter arc of $J$ between $x$ and $y$. Suppose the chord-arc property fails. By considering $C=1,2,3,\dots$ one constructs sequences $x_n,y_n \in J$ such that 
\begin{equation} \label{chordarcfailure}
\text{length}\, \big( J(x_n,y_n) \big) > n|x_n-y_n| .
\end{equation}
Notice $|x_n - y_n| \to 0$ since the length of $J(x_n,y_n)$ is bounded by the length of $J$, which is finite. Further, by compactness we may assume the sequences $x_n$ and $y_n$ converge to some point $x \in J$. The domain $\Omega$ has Lipschitz boundary by hypothesis, and so the curve $J$ can be represented near $x$ as the graph of a Lipschitz function. That is, after rotating the coordinate system suitably, there is a disk $B$ centered at $x$ and a Lipschitz function $b : \R \to \R$ such that $B \cap J = B \cap \{ (t,b(t)) : t \in \R \}$. Let $\beta$ be the Lipschitz constant. For all $n$ large enough that the points $x_n$ and $y_n$ lie in the disk $B$, choose $s_n$ and $t_n$ such that $x_n=(s_n,b(s_n))$ and $y_n=(t_n,b(t_n))$. Then 
\[
\text{length}\, \big( J(x_n,y_n) \big) \leq \sqrt{1+\beta^2}\,  |s_n-t_n| \leq \sqrt{1+\beta^2} \, |x_n-y_n| ,
\]
which contradicts \eqref{chordarcfailure} as $n \to \infty$. Therefore $J$ must satisfy the chord--arc property. 

The chord--arc property of $\partial \Omega$ implies that $\Omega$ is Ahlfors-regular \cite[Proposition 7.7]{P92}, and hence the conformal map $f$ satisfies the Smirnov condition \cite[Proposition 7.5 and Theorem 7.6]{P92}, which says that on $\D$ the harmonic function $\log |f^\prime|$ equals the Poisson integral of its boundary values. Its boundary values are constant a.e., by our work above, and so $\log |f^\prime|$ is constant on $\D$. Thus $|f^\prime|$ is constant, and so $f^\prime$ is constant, which means $f$ is linear and $\Omega$ is a disk, as we wanted to show. \qed

\bigskip
Next we justify the claim made earlier in the paper that \autoref{robinszego} fails when $\alpha>32.7$. Specifically, we show 
\begin{equation} \label{theoremfails} 
\lambda_2\big(\D;\alpha/L(\D)\big)A(\D) < \alpha \qquad \text{when $\alpha>32.7$,}
\end{equation}
so that by \eqref{asympineq} the disk $\D$ gives a smaller value than a long thin rectangle $\Omega_t$,  for large $t$, and hence the disk is not the maximizer. 

To prove \eqref{theoremfails}, recall from \autoref{preliminaries} (see \autoref{Robin_g_second}) that the second eigenfunction of the disk with positive Robin parameter $\alpha/L(\D)=\alpha/2\pi$ has radial part $g(r)=J_1(\sqrt{\lambda_2} \, r)$ where $\sqrt{\lambda_2} \in (j_{1,1}^\prime,j_{1,1})$ is chosen to satisfy the Robin boundary condition $g^\prime(1)+(\alpha/2\pi) g(1)=0$. That condition rearranges to say 
\begin{equation} \label{alphalambda} 
\alpha = -2\pi \frac{\sqrt{\lambda_2} J_1^\prime(\sqrt{\lambda_2})}{J_1(\sqrt{\lambda_2})} .
\end{equation}
Since $\lambda_2$ is a strictly increasing function of $\alpha$, we may invert and regard $\alpha$ as a function of $\lambda_2$ (see \cite[Section 5]{FL18a} with $n=2$ and $\kappa=1$). By the last formula, the condition $\alpha > \lambda_2 \pi$ in \eqref{theoremfails} is equivalent to 
\[
-2\frac{J_1^\prime(x)}{x J_1(x)} > 1 
\]
where $x=\sqrt{\lambda_2} \in (j_{1,1}^\prime,j_{1,1})$. Solving numerically, the inequality holds for $3.2261 \leq x < j_{1,1}$ where we rounded the root up to $3.2261$. Substituting this root into \eqref{alphalambda} and again rounding up, we obtain a range $32.7 \leq \alpha < \infty$ on which \eqref{theoremfails} holds.

\section{\bf Proof of \autoref{fixedalpha}}
\label{corollariesproof1}

We may assume $\Omega$ has area $\pi$, after rescaling, and so the task is to show $\lambda_2(\Omega;\alpha) \leq \lambda_2(\D;\alpha)$ when $\alpha \in [-1,0]$. 

The isoperimetric inequality $L(\Omega) \geq 2\pi$ implies $\alpha \leq 2\pi \alpha/L(\Omega)$ when $\alpha \leq 0$, and so 
\[
\lambda_2(\Omega;\alpha) \leq \lambda_2\big( \Omega; 2\pi \alpha/L(\Omega) \big) 
\]
because the Robin eigenvalues are increasing functions of $\alpha$. The assumption $\alpha \in [-1,0]$ ensures $2\pi \alpha \in [-2\pi,0]$, and so \autoref{robinszego} can be applied with $\alpha$ replaced by $2\pi \alpha$, giving  
\[
\lambda_2\big( \Omega; 2\pi \alpha/L(\Omega) \big) \leq \lambda_2\big( \D; 2\pi \alpha/L(\D) \big)  = \lambda_2(\D; \alpha) .
\]
Combining the last two inequalities proves the corollary. 

If equality holds then $\Omega$ must be a disk, by the equality statement in \autoref{robinszego}. 

\section{\bf Proof of \autoref{szegoweinstock}}
\label{corollariesproof2}

That $\mu_1(\Omega) A(\Omega)$ is maximal for the disk, under area normalization, is the case $\alpha=0$ of \autoref{robinszego}.

Weinstock's result, saying the disk maximises the first nontrivial Steklov eigenvalue under perimeter normalization, requires a little more explanation. The Steklov spectrum of the Laplacian is denoted $0=\sigma_0 < \sigma_1 \leq \sigma_2 \leq \dots$ where the eigenvalue problem is
\[
\begin{split}
\Delta u & = 0 \ \ \quad \text{in $\Omega$,} \\
\frac{\partial u}{\partial\nu} & = \sigma u \quad \text{on $\partial \Omega$.} 
\end{split}
\]
Thus $\sigma$ belongs to the Steklov spectrum exactly when $0$ belongs to the Robin spectrum with $\alpha=-\sigma$. 

After rescaling $\Omega$ we may suppose it has area $\pi$. The task is to prove $\sigma_1(\Omega) L(\Omega) \leq 2\pi$, since $\sigma_1(\D)=1$. Choosing $\alpha=-2\pi$ in \autoref{robinszego} yields
\[
\lambda_2\big( \Omega;-2\pi/L(\Omega) \big) \leq \lambda_2\big(\D;-2\pi /L(\D)\big) = \lambda_2(\D;-1) =0 .
\]
Also $\lambda_2(\Omega;0)=\mu_1(\Omega) > 0$. Since the Robin eigenvalues vary continuously with $\alpha$, a value $\widetilde{\alpha} \in [-2\pi,0)$ must exist for which $\lambda_2\big( \Omega;\widetilde{\alpha}/L(\Omega) \big) = 0$. Choose $\widetilde{\alpha}$ to be the greatest such number, so that $\lambda_2\big( \Omega;\alpha/L(\Omega) \big) > 0$ for all $\alpha > \widetilde{\alpha}$. Then $-\widetilde{\alpha}/L(\Omega)$ belongs to the Steklov spectrum of $\Omega$, and is in fact the smallest positive Steklov eigenvalue, $\sigma_1(\Omega)$. 
Hence $\sigma_1(\Omega) L(\Omega) = -\widetilde{\alpha} \leq 2\pi$, as we needed to show. 

If equality holds then $\widetilde{\alpha}=-2\pi$, and so the equality statement in \autoref{robinszego} (with $\alpha=-2\pi$) implies that $\Omega$ is a disk.

\section{\bf Proof of \autoref{mapping}}
Fix $\alpha > 0$. The Robin eigenfunction on the disk corresponding to $\lambda_1(\D;\alpha/2\pi)$ has radial form $u_1 = g(r)$, by \autoref{basic1}. Adapting P\'{o}lya and Szeg\H{o}'s method \cite{PS51}, we define  
\[
v_1 = u_1 \circ F^{-1} 
\]
on $\Omega$. This function is smooth and bounded, and belongs to $H^1(\Omega)$ by conformal invariance of the Dirichlet integral. Employing $v_1$ as a trial function in the Rayleigh principle for the first eigenvalue yields
\[
\lambda_1\big(\Omega;\alpha/L(\Omega)\big) \int_\Omega v_1^2 \, dx 
\leq \int_\Omega |\nabla v_1|^2 \, dx + \frac{\alpha}{L(\Omega)} \int_{\partial \Omega} v_1^2 \, ds , 
\]
which pulls back under the conformal map $F$ to 
\[
\lambda_1\big(\Omega;\alpha/L(\Omega)\big) \int_\D u_1^2 |F^\prime|^2 \, dx 
\leq \int_\D |\nabla u_1|^2 \, dx + \frac{\alpha}{L(\Omega)} \int_{\partial \Omega} v_1^2 \, ds ,
\]
by conformal invariance of the Dirichlet integral. 

Substituting $u_1 = g(r)$, which in particular gives $v_1=g(1)$ on $\partial \Omega$, we obtain
\begin{equation} \label{psbound1}
\lambda_1\big(\Omega;\alpha/L(\Omega)\big) \int_\D g(r)^2  |F^\prime|^2 \, dx 
\leq \int_\D g^\prime(r)^2 \, dx + \alpha g(1)^2 .
\end{equation}
For the left side of the inequality note that
\begin{equation} \label{derivativebound}
|F^\prime(0)|^2 \int_\D g(r)^2 \, r\, drd\theta \leq \int_\D g(r)^2 |F^\prime(re^{i\theta})|^2 \, r\, drd\theta  
\end{equation}
because
\begin{equation} \label{derivativebound2}
|F^\prime(0)|^2 
= \left| \frac{1}{2\pi} \int_0^{2\pi} F^\prime(re^{i\theta}) \, d\theta \right|^2 
\leq \frac{1}{2\pi} \int_0^{2\pi} |F^\prime(re^{i\theta})|^2 \, d\theta .
\end{equation}
Multiply inequality \eqref{derivativebound} by $\lambda_1\big(\Omega;\alpha/L(\Omega)\big)$, which is positive since $\alpha > 0$, and then substitute into \eqref{psbound1}, getting
\[
\lambda_1\big(\Omega;\alpha/L(\Omega)\big) \, |F^\prime(0)|^2 
\leq \frac{\int_\D g^\prime(r)^2 \, dx + (\alpha/2\pi) \int_{\partial \D} g(1)^2 \, ds}{\int_\D g(r)^2 \, dx } = \lambda_1(\D;\alpha/2\pi) ,
\]
as we wanted to prove. 

If equality holds in the theorem, then equality must hold in \eqref{derivativebound}, and hence also in \eqref{derivativebound2} for $r \in (0,1)$. By substituting the power series for $F$ into  \eqref{derivativebound2} and setting the two sides equal, we deduce that $F^\prime$ is constant and hence $F$ is linear.

\section*{Acknowledgments}
This research was supported by the Funda\c c\~{a}o para a Ci\^{e}ncia e a Tecnologia (Portugal) through project PTDC/MAT-CAL/4334/2014 (Pedro Freitas), by a grant from the Simons Foundation (\#429422 to Richard Laugesen), by travel support for Laugesen from the American Institute of Mathematics to the workshop on \emph{Steklov Eigenproblems} (April--May 2018), and support from the University of Illinois Scholars' Travel Fund. The research advanced considerably during the workshop on \emph{Eigenvalues and Inequalities} at Institut Mittag--Leffler (May 2018), organized by Rafael Benguria, Hynek Kovarik and Timo Weidl, and also benefitted from conversations with Dorin Bucur at the conference \emph{Results in Contemporary Mathematical Physics}, held in honor of Rafael Benguria in Santiago, Chile (December 2018).

\end{document}